\input amstex
\magnification=\magstep{1.5}
\baselineskip 22pt
\documentstyle{amsppt}
\topmatter
\title 
Families of hyperk\"ahler manifolds
\endtitle
\author 
Keiji Oguiso 
\endauthor
\affil 
Mathematisches Institut Universit\"at Essen, Department of Mathematical 
Sciences University of Tokyo
\endaffil
\address
D-45117 Essen Germany, 153-8914 Komaba Meguro Tokyo Japan
\endaddress
\email
mat9g0\@spi.power.uni-essen.de
\endemail
\dedicatory
Dedicated to Professor Tetsuji Shioda on the occasion of his sixtieth birthday
\enddedicatory
\abstract 
We show the density of the jumping loci of the Picard number of 
a hyperk\"ahler manifold under small deformation and provide several 
applications. In particular, we apply this to 
reveal the structure of hierarchy among all the narrow 
Mordell-Weil lattices of Jacobian K3 surfaces, which essentially reduces the 
classification of the Mordell-Weil lattices 
of Jacobina K3 surfaces - which are no more finite in number - to the 
classification of those of maximal rank $18$.     
\endabstract
\endtopmatter
\document
\head
{\S 0. Introduction - Background and Results}
\endhead 
In their article ``Families of K3 surfaces'', R. Borchards, L. Katzarkov, 
T. Pantev and N. I. Shepherd-Barron have found the following remarkable 
phenomenon concerning the behaviour of the Picard number under global 
deformation: 

\proclaim{Theorem ([BKPS])} Any complete family of minimal K\"ahler surfaces 
of Kodaira dimension 0 and constant Picard number is isotrivial. \qed
\endproclaim 

Their proof is a global argument based on the ampleness of the zero 
locus of an automorphic form on a relevant moduli space and therefore 
requires the completeness of the base space in essence.
\par
\vskip 4pt

The purpose of this notes is to generalise their Theorem to a local setting 
for a hyperk\"ahler manifold (Theorem 1) by a quite different method and 
provide several applications (Corollaries 2, 3, 6, 8 and Examples 4, 5, 7). 
\par
\vskip 4pt
A hyperk\"ahler manifold is by definition a simply connected, 
compact K\"ahlerian manifold $F$ which admits, up to scalar multiple, a unique 
non-degenerate holomorphic $2$-form $\omega_{F}$ and is one of the 
building blocks of manifolds with trivial first 
Chern class by the Bogomolov decomposition Theorem ([Be]). 
A K3 surface is nothing but a hyperk\"ahler manifold of dimension $2$. 
Due to the works of Bogomolov, Beauville and Fujiki, the following 
results which are well-known for a K3 surface also hold for a hyperk\"ahler 
manifold of any dimension: 
\roster 
\item the existence of the natural primitive integral non-degenerate 
symmetric bilinear form $(*,*)$ on $H^{2}(F, \Bbb Z)$ which induces on 
$H^{2}(F, \Bbb C) = H^{1,1}(F) \oplus \Bbb C\omega_{F} 
\oplus \Bbb C \overline{\omega}_{F}$ the polarised Hodge 
structure of weight two and satisfies $(c_{1}(L)^{2}) > 0$ if 
$L$ is ample ([Be], [Fu2]); 
\item the local Torelli Theorem for the period map given by the polarised 
Hodge structure of $H^{2}(F, \Bbb Z)$ defined in (1) ([Be]).  
\endroster
Besides original articles, we also refer the readers to [Hu1, Section 1] as 
an excellent survey about these basics on a hyperk\"ahler manifold. 
\par
\vskip 4pt

Throughout this notes, we work over the complex number field $\Bbb C$. 
We mainly consider a family of hyperk\"ahler manifolds $f : \Cal X 
\rightarrow \Delta$ over the unit disk $\Delta$ with prescribed 
``polarisation'' $\Lambda_{0}$. We regard $\Delta$ as a germ $(0 \in \Delta)$ 
of deformation of the centre fibre $F := \Cal X_{0}$ and shrink freely 
whenever it is more convenient for the statement. Therefore, the following 
three statements are equivalent to one another: 
\par
(i) $f$ is trivial as a family, that is, isomorphic to the product 
$F \times \Delta$ over $\Delta$; 
\par 
(ii) $f$ is isotrivial, that is, trivial after a finite base change 
of $\Delta$; 
\par 
(iii) all the fibres of $f$ are isomorphic.
\par
\vskip 4pt 
We denote by $\rho(F)$ the Picard number of $F$, that is, the rank of the 
N\'eron-Severi group $NS(F)$. Then $0 \leq \rho(F) \leq N := b_{2}(F) - 2$. 
For the precise definition of $\Lambda_{0}$-polarised family and for our 
purpose later, we fix an isometric marking $\tau : R^{2}f_{*} 
\Bbb Z_{\Cal X} \simeq \Lambda \times \Delta$, where $\Lambda = 
(\Lambda, (*.*))$ is a lattice of signature $(3, N - 1)$, and denote the 
resulting period map by 
$$p : \Delta \rightarrow \Cal D := \{ [\omega] \in 
\Bbb P(\Lambda \otimes \Bbb C) \vert (\omega.\omega) = 0, 
(\omega.\overline{\omega}) > 0 \} 
\subset  \Bbb P(\Lambda \otimes \Bbb C) \simeq \Bbb P^{N+1}.$$ 
This map is defined by $p(t) := \tau (\omega_{\Cal X_{t}})$ and is known to 
be holomorphic. In the case of K3 surfaces, $\Lambda$ is the K3 
lattice $\Lambda_{\text{K3}} := U^{\oplus 3} 
\oplus E_{8}(-1)^{\oplus 2}$ and we have $N = 20$ and $0 \leq \rho(F) 
\leq 20$. We call a family $\Lambda_{0}$-polarised if $\tau^{-1}
(\Lambda_{0}) \subset NS(\Cal X_{t})$ for all $t \in \Delta$. We do not 
require hyperbolicity of $\Lambda_{0}$ in apriori and even allow 
$\Lambda_{0}$ to be $\{0\}$, the later of which is nothing but the case of no 
fixed "polarisations". By the term $I$-topology, we mean the topology on 
$\Delta$ for which the collection of open sets is
$\{\, \emptyset, \, \Delta, \, \Delta - \{\text{countably\, many\,  
points}\} \}$. The $I$-toplogy is compatible with shrinking. We call a real 
valued function $y = a(x)$ on a 
topological space $X$ strictly upper semi-continuous at $x_{0} \in X$ 
if there exists an open neighborhood $U$ of $x_{0}$ such that $a(x_{0}) > 
a(x')$ for all $x' \in U - \{x_{0}\}$. 
\par
\vskip 4pt

Our main Theorem is as follows: 

\proclaim{Theorem 1} For any non-trivial, $\Lambda_{0}$-polarised family 
$f : \Cal X \rightarrow \Delta$ of hyperk\"ahler manifolds, there exists a 
dense, countable subset $\Cal S \subset \Delta$ in the classical topology 
such that the function $\rho(t) := \rho(\Cal X_{t})$ on $\Delta$ is strictly 
upper semi-continous at each point of $\Cal S$ in the $I$-topology but is 
constant on the complement $\Delta - \Cal S$, provided that all the fibres 
are algebraic in the case of $\dim F \geq 4$. 
\endproclaim 

It is an easy fact that $\Cal S$ is at most countable and the essence 
of Theorem 1 is in the converse: the existence of enough jumping points.
Note that by Baire's category Theorem, $\Delta - \Cal S$ is 
also dense but uncountable and is "much bigger" than $\Cal S$. 
Note also that the same problem does not make much sense for 
Calabi-Yau manifolds: Indeed, we have $\text{Pic}(X) = H^{2}(X, \Bbb Z)$ 
for them.
\par
\vskip 4pt  
Our proof is based on the materials (1) and (2) above and is quite primitive. 
For this reason, most part of our proof 
goes through regardless of the dimension of fibres. The only difference at the 
present between the cases of dimension 2 and higher is a regrettable lack 
of the criterion for algebraicity in higher dimension: Whether or not the 
hyperbolicity of $NS(F)$ implies the algebraicity of $F$. If this is 
affirmative, then the assumption on algebraicity made at the end of Theorem 1 
is not needed. Refer to [Hu2] for positive directions. We shall prove 
Theorem 1 in Section 1.
\par
\vskip 4pt 
Thanks to the existence of the coarse moduli scheme for polarised manifolds 
due to Viehweg ([Vi]), Theorem 1 deduces the following slightly 
stronger isotriviality in an algebraic setting:

\proclaim{Corollary 2} Any smooth family $g : \Cal Y \rightarrow \Cal B$ 
projective over a normal noetherian irreducible base scheme $\Cal B$ of 
hyperk\"ahler manifolds with constant Picard number becomes trivial after 
appropriate \'etale Galois covering $\Cal B' \rightarrow \Cal B$ of the base 
scheme. In particular, if in addition $\pi_{1}^{\text{alg}}(\Cal B) = \{1\}$, 
then $g$ itself is trivial.  
\endproclaim

Here we do not require the completeness of $\Cal B$. 
\par 
\vskip 4pt 
It is known that there exists a non-isotrivial, smooth projective family of 
supersingular K3 surfaces over $\Bbb P^{1}$. Therefore the corresponding 
statement in positive characteristic is false even if we assume the base 
space to be complete. (Refer to [GK] for relevant phenomena in positive 
characteristics.) We shall prove Corollary 2 in Section 2. 
\par
\vskip 4pt
Besides geometrical applications, Theorem 1 also deduces the following 
lovely result on arithmetic. We shall prove this in Section 2:

\proclaim{Corollary 3} Let $\Bbb H$ be the upper half plane and 
$\text{GL}^{+}(2, \Bbb Q)$ the index two 
subgroup of $\text{GL}(2, \Bbb Q)$ consisting of elements $M$ such that 
$\text{det}M > 0$. Let $w = \varphi(z)$ be a holomorphic 
function defined over a neighbourhood of $\tau \in \Bbb H$.
Assume that $\varphi(\tau) \in \Bbb H$. Then, there exists a sequence 
$\{\tau_{k}\}_{k=0}^{\infty} \subset {\Bbb H} - \{\tau\}$ such that 
$\lim_{k \rightarrow \infty} \tau_{k} = \tau$ and that $\tau_{k}$ and 
$\varphi(\tau_{k})$ are congruent for each $k$ with respect to the standard,  
linear fractional action of $\text{GL}^{+}(2, \Bbb Q)$ on $\Bbb H$.  
\endproclaim  
In the rest of this Introduction, we shall restrict ourselves to a family 
of K3 surfaces and see what happens more closely. Explicit construction and 
proof will be given in Section 3 and Section 2 respectively. 
\par
\vskip 4pt
It is reasonable to measure how algebraic a given K3 surface is by the 
lexicographic order of pairs $(a(S), \rho(S))$. Here $a(S)$ is the algebraic 
dimension. From this view point, the next example might be of some interest:

\proclaim{Example 4} There is a family of K3 surfaces 
$f : \Cal X \rightarrow \Delta$ for which there exists a dense, countable 
subset $\Cal S \subset \Delta$ such that $a(\Cal X_{s}) = 2$ and 
$\rho(\Cal X_{s}) = 20$, the maximum, for $s \in \Cal S$, but 
$a(\Cal X_{t}) = 0$, the minimum, and $\rho(\Cal X_{t}) = 19$ 
for $t \not\in \Cal S$.  
\endproclaim

The values of $\rho$ do not always sweep out some range of integers 
``continuously'': 

\proclaim{Example 5} There is a family of algebraic K3 surfaces $f : \Cal X 
\rightarrow \Delta$ for which there exists a dense countable subset $\Cal S$ 
such that $\rho(\Cal X_{s}) = 20$ for $s \in \Cal S$, 
but $\rho(\Cal X_{t}) = 18$ for $t \not\in \Cal S$. In particular, there 
lacks the intermediate integer $19$. 
\endproclaim
\par
\vskip 4pt
A Jacobian K3 surface is a K3 surface equipped with an elliptic fibre space 
structure $\varphi : S \rightarrow \Bbb P^{1}$ with section $O$. By a local 
family of Jacobian K3 surfaces, we mean a commutative diagram 
$$\CD 
\Cal X @> \varphi >> \Cal W @> \Cal O >> \Cal X \\
@V f VV @VV \pi V @VV f V \\
\Delta @> id >> \Delta @> id >> \Delta
\endCD$$ 
such that $\varphi_{t} : \Cal X_{t} \rightarrow \Cal W_{t}$ is a Jacobian 
K3 surface with section $\Cal O_{t}$. 
\par
\vskip 4pt
The Mordell-Weil lattice $MW(\varphi)$ of $\varphi$ is the Mordell-Weil 
group, that is, the group of sections of $\varphi$, equipped with Shioda's 
positive definite, symmetric bilinear form 
$\langle *, *\rangle$ [Sh2]. This lattice 
structure on $MW(\varphi)$ made the study of Mordell-Weil groups extremally 
rich [Sh3]. We refer the readers to [Sh2] and [Sh3] for Mordell-Weil lattices. 
We denote by $r(\varphi)$ the rank of $MW(\varphi)$. Note that 
$0 \leq r(\varphi) \leq 18$ for a Jacobian K3 surface [Sh1].
\par
\vskip 4pt 
A similar but slightly different jumpinig phenomenon is found for 
$r(\varphi)$: 

\proclaim{Corollary 6} For any family $f = \pi \circ \varphi: 
\Cal X \rightarrow \Cal W \rightarrow \Delta$ of Jacobian K3 surfaces,  
there exists a dense countable subset $\Cal S' \subset \Delta$ in the 
classical topology such that the function $r(t) := r(\Cal X_{t})$ is strictly 
upper semi-continous at each $s \in \Cal S'$ in the $I$-topology unless 
$f : \Cal X \rightarrow \Delta$ is trivial. 
\endproclaim

\proclaim{Example 7} There exists a family $f = \pi \circ \varphi: \Cal X 
\rightarrow \Cal W \rightarrow \Delta$ of Jacobian K3 surfaces such that 
$\rho(0) = 20$ but $r(0) = 0$, and that $\rho(t) < 20$ but $r(t) > 0$ for 
general $t$ in the sense of the $I$-topology. In particular, $\rho(t)$ is 
strictly upper semi-continous at 
$t = 0$ but $r(t)$ is strictly lower semi-continuous at $t = 0$ in the 
$I$-topology. 
\endproclaim 
Example 7 also shows that there is a case where the behaviour of $r(t)$ is not 
honestly accompanied with that of $\rho(t)$. For comparison, we also note 
that $\rho(t)$ and $r(t)$ for a family of rational Jacobian surfaces are 
constant and lower semi-continuous respectively by the stability Theorem. 
\par
\vskip 4pt
Our final aim of this notes is to apply Theorem 1 to study the structure of 
the Mordell-Weil lattices of Jacobian K3 surfaces. By the narrow Mordell-Weil 
lattice $MW^{0}(\varphi)$ we mean the sublattice of $MW(\varphi)$ of finite 
index consisting of the sections which pass through the identity component of 
each fibre [Sh2]. Contrary to the case of rational Jacobian 
surfaces, the isomorphism classes of both $MW(\varphi)$ and $MW^{0}(\varphi)$ 
for Jacobian K3 surfaces are no more finite ([OS], [Ni]) and the whole 
pictures of them does not seem so clear even now. Our interest here is 
to clarify certain relationships among all of $MW(\varphi)$ through Theorem 1:

\proclaim{Corollary 8} For any given Jacobian K3 surface $\varphi : J 
\rightarrow \Bbb P^{1}$ of rank $r := r(\varphi)$, there exists a sequence 
$\{\varphi_{m} : J_{m} \rightarrow \Bbb P^{1}\}_{m = r}^{18}$ of Jacobian K3 
surfaces such that 
\roster
\item 
$\varphi_{r} : J_{r} \rightarrow \Bbb P^{1}$ is 
the original $\varphi : J \rightarrow \Bbb P^{1}$;
\item 
$r(\varphi_{m}) = m$ for each $m$; and 
\item 
there exists a sequence of isometric embeddings:
\endroster
$$MW^{0}(\varphi) (= MW^{0}(\varphi_{r})) \subset MW^{0}(\varphi_{r+1}) 
\subset \cdots \subset MW^{0}(\varphi_{17}) \subset MW^{0}(\varphi_{18}).$$   
In particular, the narrow Mordell-Weil lattice of any Jacobian K3 surface is 
embedded into the Mordell-Weil lattice of some Jacobian K3 
surface of rank $18$. Conversely, for every sublattice $M$ of the narrow 
Mordell-Weil lattice of a Jacobian K3 surface of rank $18$, there exists 
a Jacobian K3 surface whose narrow Mordell-Weil lattice contains 
$M$ as a sublattice of finite index. Moreover, for each given $M$ 
there are at most finitely many isomorphism classes of the Mordell-Weil 
lattices of Jacobian K3 surfaces which contains $M$ as a sublattice of 
finite index.  
\endproclaim 
This Corollary coarsely reduces the study of $MW(\varphi)$ to those of 
the maximal rank $18$ and is in some sense an analogue of the fact 
that $MW^{0}(\varphi)$ of rational Jacobian surfaces are embedded into 
$E_{8}$. Concerning rational Jacobian case, see [OS] together with a 
correction at the end of Section 3 for the classification, and 
refer to [Sh3, Part II] for beautiful aspects behind the hierarchy governed 
by the lattice $E_{8}$. It might be also worthwhile noticing here that a 
Jacobian K3 
surface with Mordell-Weil rank 18 is necessarily ``singular'' 
in the sense of Shioda [SI] and that Nishiyama [Ni] has already constructed 
an infinite series of examples of such Jacobian K3 surfaces. 
We shall prove Corollary 8 in Section 2.  
\par
\vskip 4pt 
\subhead
Acknowledgement 
\endsubhead
\par 
\vskip 4pt
Besides pioneering works [BKPS] and [Sh2], a greater part of idea of this 
notes has been brought to the author by Professor Eckart Viehweg through their 
discussion concerning the Shafarevich Conjecture (in progress) and the idea 
of Example 4 has grown out through the author's joint work with Professor 
Thomas Peternell [OP]. First of all, the author would like to express his 
deep thanks to both of them for their extremally valuable conversations, 
without which the author could never carry out this work. He also would like 
to express his thanks to T. Szemberg, E. Bedoulev for their valuable 
discussion. Especially, some question of T. Szemberg to the author brought to 
him Corollary 3 and Example 5. This work has been done during the author's 
stay in the 
Universit\"at Essen 1999 under the financial support by the 
Alexander-Humboldt fundation. It is their generous support of Professors 
H. Esnault, F. Catanese, Y. Kawamata, D. Morrison, V. V. Shokurov and E. 
Viehweg that made the author's stay in Germany in the academic year 1999 
possible. The author would like to express his sincere thanks to all of them. 
The author would like to express his sincere thanks to the Alexander-Humboldt 
fundation for the financial support, and Universit\"at Essen for providing a 
pleasant environment to research. Taking this 
place, the author would like to express his deep thanks to Professor T. 
Katsura for informing him a countre-example in positive characteristic 
and to Professor Y. Yamada for not only pointing out mistakes in [OS] but 
also showing him correction of them, which the author quoted at the end of 
Section 3. Last but not least at all, the author would like to express his 
hearty thanks to Professors Yujiro Kawamata, Hajime Kaji and also Miss Kaori 
Suzuki 
for their warm encouragement through e-mails. It is Professor Tetsuji Shioda 
who introduced to the author 
his beautiful theory of Mordell-Weil lattices, on which this notes depends 
much. It would be a great pleasure for the author to dedicate this notes to 
Professor Tetsuji Shioda on the occasion of his sixtieth birthday.   
\par 
\vskip 4pt 
\head
{\S 1. Proof of Theorem 1}
\endhead 
We employ the same notation and convention as in Introduction. 
\proclaim{Lemma (1.1)}
There exists a unique primitive sublattice $\Lambda_{1} \subset 
\Lambda$ such that: 
\roster 
\item 
$\Lambda_{0} \subset \Lambda_{1} \subset \tau(NS(\Cal X_{t}))$ for all $t \in 
\Delta$, and that;
\item 
if $\Lambda \subset \tau(NS(\Cal X_{t}))$ for all $t \in \Delta$, then 
$\Lambda \subset \Lambda_{1}$. 
\endroster
Moreover, there exists a unique, non-empty $I$-open subset $\Cal U$ such that 
$\tau(NS(\Cal X_{t})) = \Lambda_{1}$ for ecah $t \in \Cal U$ and that 
$\tau(NS(\Cal X_{t}))$ is strictly bigger than $\Lambda_{1}$ if $t \not\in 
\Cal U$. In particular, the jumping points of $\rho$ is at most 
countable.   
\endproclaim 
\demo{Proof} Consider all the primitive sublattices $\Lambda_{n}$ ($n \in 
\Cal N$) of $\Lambda$ which contain $\Lambda_{0}$ and set $\Delta(n) := 
\{t \in \Delta \vert \tau(NS(\Cal X_{t})) = \Lambda_{n} \}$. Then $\Delta = 
\sqcup_{n \in \Cal N} \Delta(n)$, because $NS(\Cal X_{t})$ is primitive in 
$H^{2}(\Cal X_{t}, \Bbb Z)$. Recall that $\Cal N$ is countable but 
$\Delta$ is uncountable and that countable union of countable sets is again 
countable. Therefore there exists an element of $\Cal N$, which we denote by 
$1$, such that $\Delta(1)$ is uncountable. By the Lefschetz (1,1)-Theorem, 
$p(t) \in \Lambda_{1}^{\perp} \otimes \Bbb C$ for all $t \in \Delta(1)$. Then 
$p(\Delta)$ must be contained in the linear space defined by $(\Lambda_{1}.*) 
= 0$ in $\Bbb P(\Lambda \otimes \Bbb C)$, because $p$ is holomorphic. 
Therefore, again, by the Lefschetz (1,1)-Theorem, we have $\Lambda_{1} 
\subset \tau(NS(\Cal X_{t}))$ for all $t \in \Delta$. Set $\Cal U := \{t \in  
\Delta \vert \tau(NS(\Cal X_{t})) = \Lambda_{1} \}$. Then, by the choice of 
$\Lambda_{1}$, we know that $\Cal U \not= \emptyset$. In addition, since 
$\Lambda_{1} \subset \tau(NS(\Cal X_{t}))$ for all $t \in \Delta$ and since 
$\Lambda_{1}$ is primitive, we see that  $\tau(NS(\Cal X_{t})) \not= 
\Lambda_{1}$ if and only if there exists an element $v \in 
\Lambda - \Lambda_{1}$ such that $p(t) \in v^{\perp}$, that is, $(p(t).v) = 
0$. Set $\Cal S := 
\cup_{v \in \Lambda - \Lambda_{1}}(p(\Delta) \cap v^{\perp})$. 
Then $\Cal U = \Delta - p^{-1}(\Cal S)$. Since $\Cal U \not= \emptyset$, we 
have $p(\Delta) \not= p(\Delta) \cap v^{\perp}$ for each $v$. Since $p$ is 
holomorphic, this implies that $p(\Delta) \cap v^{\perp}$ is at 
most countable. Hence $\Cal S$ is also countable. Therefore, $\Cal U$ is 
$I$-open, again by the fact that $p$ is holomorphic. Now we are 
done. \qed
\enddemo 
As in Introduction, we set $N := b_{2}(F) - 2$, where $F$ is the centre fibre 
of $f$. By abuse of notation, we denote the rank of a lattice $L$ by 
$\dim L$. 
\proclaim{Lemma (1.2)} Let $\Lambda_{1}$ be the lattice found in (1.1). 
Then $\dim \Lambda_{1} \leq N - 1$. In other words, $f$ is trivial if 
$\dim \Lambda_{1} = N$. 
\endproclaim    
\demo{Proof} We have $\dim \Lambda_{1} \leq N$. Assume 
for a contradiction that $\dim \Lambda_{1} = N$. Then $\Lambda_{1} 
\otimes \Bbb R \simeq H^{1, 1}(F, \Bbb R)$ and $\Lambda_{1}^{\perp} \otimes 
\Bbb R \simeq \Bbb R \langle \text{Re}(\omega_{F}), \text{Im}(\omega_{F}) 
\rangle$ by the dimension reason. In particular, 
$\Lambda_{1}^{\perp}$ is positive definite and is of rank $2$. In addition, 
since $(p(t).\Lambda_{1}) = 0$, we have $p(\Delta) \subset  \{ [\omega] \in 
\Bbb P(\Lambda_{1}^{\perp} \otimes \Bbb C) \vert (\omega.\omega) = 0, 
(\omega.\overline{\omega}) > 0 \} \subset  \Bbb P(\Lambda_{1}^{\perp} \otimes 
\Bbb C) \simeq \Bbb P^{1}$, where we regard $\Bbb P(\Lambda_{1}^{\perp} 
\otimes \Bbb C)$ as the linear subspace of $\Bbb P(\Lambda \otimes \Bbb C)$ 
defined by $(\Lambda_{1}.*) = 0$. Note that the equation 
$(\omega.\omega) = 0$ has at most 
two solutions in $\Bbb P^{1}$,  because $(*,*)$ is positive definite on 
$\Lambda_{1}^{\perp}$ and therefore is non-degenerate. Since $\Delta$ is 
connected and $p$ is continuous, this implies that $p(\Delta)$ consists of 
one point. That is, the period map $p$ is constant. Hence, by the local 
Torelli Theorem, the fibres of $f$ are all isomorphic and $f$ is then trivial 
around $0$. However this contradicts our assumption. \qed
\enddemo 
Set $\Cal S := \Delta - \Cal U$, where $\Cal U$ is same as in (1.1). 
The next Lemma completes the proof: 
\proclaim{Lemma (1.3)} $\Cal S$ is dense in $\Delta$ in the classical 
topology.  
\endproclaim 
\definition{Convention for argument of (1.3) and (1.4)} 
In order to prove (P) by the argument by contradiction, we 
first assume that (P) does not hold and proceed our argument until we reach 
a contradiction. However, we often meet a situation where we need to prove 
another property (Q) again by the argument by contradiction. 
Assume (Q) does not hold, then $\cdots$. 
In this case, if we reach a contradiction, then the logical conclusion at 
this stage should be written: Either (P) or (Q) does hold. However, if (P) 
does hold but (Q) does not hold, then this itself gives a contradiction to 
the logic we have proceeded, because this time ``(Q) does not hold'' is the 
truth, and we have nothing to do more! For this reason, and 
especially for simplicity of description, our proof 
below is written as if the conclusion at that stage would be that: 
the {\it last} property (Q) does hold. Please accept this logically 
harmless convention in the proof below. 
\enddefinition  
\demo{Proof} Assume to the contrary that $\Cal S$ is not dense. Then, 
$\overline{\Cal S} \not= \Delta$, where the closure is taken inside $\Delta$. 
Since $\overline{\Cal S} \subset \Delta$, there then exists a point 
$P \in \Delta - \overline{\Cal S}$. Since $\Delta - \overline{\Cal S}$ is 
open in $\Delta$, by the definition of the induced topology, there exists an 
open subset $U$ of $\Bbb C$ such that $U \cap \Delta = \Delta - 
\overline{\Cal S}$. Since $\Delta$ is open in $\Bbb C$, we see that $U \cap 
\Delta$ is also open in $\Bbb C$. Therefore, there exists a small disk 
$\Delta_{P}$ centered at $P$ such that $P \in  \Delta_{P} \subset \Delta - 
\overline{\Cal S}$. Then, by the definition of $\Cal S$ and by the inclusions 
$\Delta_{P} \subset \Delta - \overline{\Cal S} \subset \Delta - \Cal S = 
\Cal U$, we see that $\tau(NS(\Cal X_{t})) = \Lambda_{1}$ for all $t \in 
\Delta_{P}$ by (1.1). 
\enddemo 
\proclaim{Claim (1.4)} $p \vert \Delta_{P}$ is constant.
\endproclaim 
Once we get (1.4), then the result (1.3) follows. Indeed, constantness of $p$ 
on $\Delta_{P}$ implies the constantness of $p$ on the whole $\Delta$, because 
$p$ is holomorphic. Then, by the local Torelli Theorem, all the fibres are 
isomorphic. This in particular implies that $p$ is trivial around $0$, a 
contradiction. 

\demo{Proof of (1.4)} Assume to the contrary that $p \vert \Delta_{P}$ is not 
constant. Assume that $\Cal X_{P'}$ of dimension 2 is not algebraic for 
some $P' \in \Delta_{P}$. Then there exists $P'' \in \Delta_{P}$ such that 
$\Cal X_{P''}$ is algebraic. This is due to [Fu1, Theorem 4.8 (2)] and is 
valid in any dimension. Then on the one hand, $\Lambda_{1} \simeq 
NS(\Cal X_{P'})$ is negative semi-definite by the algebraicity criterion 
for surface, but on the other hand $\Lambda_{1} \simeq NS(\Cal X_{P''})$ is 
hyperbolic, a contradiction. Therefore, $\Cal X_{P'}$ is algebraic for all 
$P' \in \Delta_{P}$ even in the case of dimension 2. 
{\it From now on, dimension of fibre is arbtrary.} Note that a K\"ahler 
Moishezon manifold is projective. Therefore, by the algebraicity of 
fibre, $\Lambda_{1}$ contains an element $a$ such that $(a, a) > 0$. 
Then by the Schwartz inequality, $\Lambda_{1}$ is non-degenerate and even 
hyperbolic, because $\Lambda_{1} \subset \tau(H^{1, 1}(F, \Bbb R))$ and 
$\tau(H^{1, 1}(F, \Bbb R))$ is of index $(1, *)$. Let us choose a 
holomorphic coordinate $z$ of $\Delta_{P}$ centered at $P$. By the Lefschetz 
(1,1)-Theorem, we have $p \vert \Delta_{P} : \Delta_{P} \rightarrow  
\{ [\omega] \in 
\Bbb P(\Lambda_{1}^{\perp} \otimes \Bbb C) \vert (\omega.\omega) = 0, 
(\omega.\overline{\omega}) > 0 \} \subset  \Bbb P(\Lambda_{1}^{\perp} \otimes 
\Bbb C) \simeq \Bbb P^{n}$, where we again regard $\Bbb P(\Lambda_{1}^{\perp} 
\otimes \Bbb C)$ as a linear subspace of $\Bbb P(\Lambda \otimes \Bbb C)$ 
defined by $(\Lambda_{1}.*) = 0$. Note that $n \geq 2$ by (1.2). Let us fix an 
integral  basis of $\Lambda_{1}^{\perp}$ and write $p\vert \Delta_{P}$ as 
$p(z) = [F_{0}(z) : F_{1}(z) : F_{2}(z) : \cdots : F_{n}(z)]$. We may assume 
without loss of generality that $F_{0}(0) \not= 0$. Then we may rewrite $p$ 
under the inhomogeneous coordinates $x_{i} := X_{i}/X_{0}$ as $z \mapsto 
[1 : f_{1}(z) :  f_{2}(z) : \cdots : f_{n}(z)]$ around $P$. In this 
expression, $f_{i}(z)$ are holomorphic functions of $z$. 
Assume that $f_{i}(z)$ are all constant. Then $p$ is constant around $P$ 
and therefore constant on the whole $\Delta_{P}$. However, this contradicts 
our assumption that $p$ is not constant on $\Delta_{P}$. Therefore, some 
$f_{k}(z)$ is not constant. Then by changing the order of the coordinates, 
we may assume without loss of generality that $f_{1}(z)$ is not 
constant. Since $\text{dim}_{\Bbb R}\Bbb C 
= 2$, the three elements $1, f_{1}(0), f_{k}(0)$ in $\Bbb C$ are linearly 
dependent over $\Bbb R$ for each $k \in \{2, \cdots, n\}$. Therefore, there 
exists $(r_{0,k}, r_{1,k}, r_{2,k}) 
\in \Bbb R^{3} - \{\overrightarrow{0}\}$ such that $r_{0,k}\cdot1 + 
r_{1,k}f_{1}(0) + r_{2,k}f_{k}(0) = 0 - (*)$. In what follows, for 
$\overrightarrow{a} := 
(a_{0}, a_{1}, a_{2}, 
\cdots, a_{n}) \in \Bbb R^{n+1} - \{\overrightarrow{0}\}$, 
we denote the hyperplane defined by 
$a_{0}X_{0} + a_{1}X_{1} + a_{2}X_{2} + \cdots + a_{n}X_{n} = 0$ 
in $\Bbb P(\Lambda_{1}^{\perp} \otimes \Bbb C)$ by $H(\overrightarrow{a})$ 
and put $f_{\overrightarrow{a}}(z) := 
a_{0} + a_{1}f_{1}(z) + a_{2}f_{2}(z) + \cdots + a_{n}f_{n}(z)$. 
Set $\overrightarrow{r_{k}} := (r_{0,k}, r_{1,k}, 0, \cdots, 0, r_{2,k}, 0, 
\cdots, 0)$, where $r_{2,k}$ is located at same position as $X_{k}$ 
in $[X_{0} : X_{1} : \cdots : X_{n}]$. 
In this terminology, the previous equation (*) is equivalent 
to $p(0) \in H(\overrightarrow{r_{k}})$ and is also equivalent to 
$f_{\overrightarrow{r_{k}}}(0) = 0$. In particular, $p(\Delta_{P}) \cap 
H(\overrightarrow{r_{k}}) \not= \emptyset$. 
Assume that for some $k$, we have $p(\Delta_{P}) \not\subset 
H(\overrightarrow{r_{k}})$, that is, $f_{\overrightarrow{r_{k}}}(z) 
\not\equiv 0$. Since $0$ is a zero of $f_{\overrightarrow{r_{k}}}(z)$, by 
using the discreteness of zeros of a non-zero holomorphic 
function, we may choose a small circle $\gamma \subset \Delta_{P}$ around 
$z = 0$ such that $f_{\overrightarrow{r_{k}}}(z)$ has no zeros on $\gamma$. 
Set $K := \min \{\vert f_{\overrightarrow{r_{k}}}(z) \vert \vert z \in 
\gamma\}$ and 
$M := \max \{\vert f_{i}(z) \vert z \in \gamma, i = 0, 1, \cdots, n\}$, 
where we define $f_{0}(z) \equiv 1$. Note that $K > 0$ and $M > 0$. 
Then, by using the triangle inequality, we see that 
$\vert f_{\overrightarrow{r_{k}}}(z) - f_{\overrightarrow{a}}(z) \vert 
< \vert f_{\overrightarrow{r_{k}}}(z) \vert$ on $\gamma$ 
provided that $\vert \overrightarrow{a} - \overrightarrow{r_{k}} \vert < 
KM^{-1}(n+1)^{-1}$. Denote by $U$ the open disk such that $\partial U = 
\gamma$. Then, by Rouch\'e's Theorem, the cardinalities 
of zeros counted with multiplicities on $U$ 
are the same for $f_{\overrightarrow{r_{k}}}$ and 
$f_{\overrightarrow{a}}$. Therefore, 
$f_{\overrightarrow{a}}$ also admits a zero on $U$, 
because $0$ is a zero of $f_{\overrightarrow{r_{k}}}$.  
Since $\Bbb Q^{n+1} - \{\overrightarrow{0}\}$ is dense in 
$\Bbb R^{n+1} - \{\overrightarrow{0}\}$ and since $KM^{-1}(n+1)^{-1} > 0$, 
there then exists $\overrightarrow{q} \in \Bbb Q^{n+1} - 
\{\overrightarrow{0}\}$ such that $f_{\overrightarrow{q}}(z)$ admits 
a zero on $U$. Let us denote this zero by $Q \in U (\subset \Delta_{P})$. 
Then, $f_{\overrightarrow{q}}(Q) = 0$. Or in other words, 
$p(Q) \in H(\overrightarrow{q})$.  
Recall that $\Lambda_{1}^{\perp}$ is 
non-degenerate because so is $\Lambda_{1}$. Note also that our homogeneous 
coordinates $[X_{0} : X_{1} : \cdots : X_{n}]$ is chosen with respect to an 
integral basis of $\Lambda$ and the rational linear equations 
$(\Lambda_{1}.*) = 0$. Therefore, there exists a rational 
vector $0 \not = v \in \Lambda_{1}^{\perp} \otimes \Bbb Q$ such that $(v.*)$ 
gives the same linear space as $H(\overrightarrow{q})$. Then 
$(v.p(Q)) = 0$. This together with the Lefschetz (1,1)-Theorem implies that 
$mv \in \tau(NS(\Cal X_{Q}))$, where $m > 0$ is a suitable integer. Note also 
that $mv \not\in \Lambda_{1}$ again because of the non-degeneracy of 
$\Lambda_{1}$. Then $\tau(NS(\Cal X_{Q}))$ becomes strictly bigger than 
$\Lambda_{1}$. On the other hand, since $Q \in \Delta_{P} \subset 
\Delta - \Cal S = \Cal U$, we have $\tau(NS(\Cal X_{Q})) = \Lambda_{1}$ 
by (1.1), a contradiction. Therefore $p(\Delta_{P})$ 
is contained in the intersection of the 
$(n-1)$-hyperplanes $H(\overrightarrow{r_{k}})$ given by 
$r_{0, k}X_{0} + r_{1, k}X_{1} + r_{2, k}X_{k} = 0$, 
$k = 2, 3, \cdots n$. In other words, we have  
$r_{0, k} + r_{1, k}f_{1}(z) + r_{2, k}f_{k}(z) \equiv 0$ for all $k$. 
Assume that $r_{2, k} = 0$ for some $k$. Then $r_{0, k} + r_{1, k}f_{1}(z) 
\equiv 0$. 
Therefore $r_{0,k} = r_{1,k} = 0$ by the non-constantness of $f_{1}(z)$. 
However, this contradicts $(r_{0,k}, r_{1,k}, r_{2,k}) \not = (0, 0, 0)$. 
Therefore $r_{2, k} \not = 0$ for any $k$. 
Hence the intersection of these $(n-1)$-hyperplanes 
$H(\overrightarrow{r_{k}})$ is a line $L \simeq \Bbb P^{1}$ by the rank 
reason. Note that this $L$ is defined over $\Bbb R$ in $\Bbb P^{n}$. Then, 
we may choose a real basis of $\Lambda_{1}^{\perp} 
\otimes \Bbb R$ under which the coordinate description of $p$ is of the form 
$p(z) = [1 : g(z) : 0 : \cdots : 0]$. Assume that $g(z)$ is constant. 
Then $p$ is also constant, a contradiction to the assumption that $p(z)$ is 
constant on $\Delta_{P}$. Therefore $g(z)$ is not constant. 
 Let us consider the matrix representation of 
the symmetric bilinear form $(*.*)\vert\Lambda_{1}^{\perp}$ under this real 
basis and denote it by $(x, y) = x(c_{ij})x^{t}$. Here $(c_{ij})$ is a real 
symmetric matrix. Since $(p(z).p(z)) = 0$, we have 
$$c_{00} + 2c_{10}g(z) + c_{11}g(z)^{2} \equiv 0 \, (**).$$
Differentiating (**) by $z$, we get $(c_{10} + c_{11}g(z))g'(z) \equiv 0$. If 
$g'(z) \equiv 0$, then $g(z)$ is constant, a contradiction. Hence $c_{10} + 
c_{11}g(z) \equiv 0$, because $\Bbb C\{z\}$ is an integral domain. Since 
$g(z)$ is not constant, this implies that $c_{11} = c_{10} = 0$. 
Then, we also get $c_{00} = 0$ by substituting $c_{11} = c_{10} = 0$ to 
the equation (**). Therefore, we have $(p(z).\overline{p(z)}) = c_{00} + 
c_{10}(g(z) + \overline{g(z)}) + c_{11}g(z)\overline{g(z)} = 0 + 0 + 0 = 0$. 
Since we changed coordinates only over 
$\Bbb R$, we still have 
$[\overline{\omega_{\Cal X_{z}}}] = [\overline{p(z)}] = [1 : 
\overline{g_{1}}(z) : 0 : \cdots 0]$. 
Therefore, we have $(p(z).\overline{p(z)}) = 
\vert c(z) \vert^{2}(\omega_{\Cal X_{z}}. \overline{\omega_{\Cal X_{z}}})$, 
where $c(z)$ is a nowhere vanishing complex function adjusting ambiguities of 
the choice of ratio. Since $(\omega_{\Cal X_{z}}. 
\overline{\omega_{\Cal X_{z}}}) > 0$ by the definition of the period map, 
we have then $(p(z).\overline{p(z)}) > 0$. However, this contradicts 
the previous equality $(p(z).\overline{p(z)}) = 0$. Now we are done. \qed 
\enddemo
\remark{Remark (1.5)} In general, given a non-constant holomorphic map 
$g : \Delta 
\rightarrow \Bbb P^{n}$, the condition $h(\Delta) \cap H(\overrightarrow{a}) 
\not = \emptyset$ is not open around $\overrightarrow{a}$ 
if this point satisfies $h(\Delta) \subset H(\overrightarrow{a})$. For 
example, consider the holomorphic map 
$h(z) = [1 : -1 : z]$. Then $h(\Delta) \subset H((1, 1, 0))$ but $h(\Delta) 
\cap H((1-\epsilon, 1, 0)) = \emptyset$ for any $\epsilon \not= 0$. 
\qed
\endremark 
\par 
\vskip 4pt 
\head
{\S 2. Proof of Corollaries}
\endhead 
In this section, we shall prove all Corollaries mentioned in Introduction. 
Here we again employ the same notation and convention as in Introduction. 

\demo{Proof of Corollary 2} Let us first show that all the fibres of $g$ are 
isomorphic. Only for this purpose, we may assume $\Cal B$ is smooth because 
of the existence of resolution due to Hironaka and may work in the complex 
analytic category. Fix $P \in \Cal B$ and take an arbitrary $Q \in \Cal B$. 
Then we may join $P$ and $Q$ by a chain of small disks $\Delta_{i}$, because 
$\Cal B$ is now assumed to be smooth and is noetherian and connected. By 
applying the contraposition of Theorem 1 for each 
$g \vert g^{-1}(\Delta_{i}) : g^{-1}(\Delta_{i}) \rightarrow 
\Delta_{i}$ and using the assumption of constantness of Picard numbers, 
we see that all the fibres of $g \vert g^{-1}(\Delta_{i})$ are isomorphic 
for each $\Delta_{i}$. Therefore $\Cal Y_{Q}$ 
is isomorphic to $\Cal Y_{P}$ as a complex manifold. 
Since both are projective, $\Cal Y_{Q}$ is then isomorphic 
to $\Cal Y_{P}$ as an algebraic variety by Chow's Theorem. {\it 
From now on, we stick to the original given base scheme} 
$\Cal B$ {\it and work in the algebraic category.} Let us denote by $F$ the 
hyperk\"ahler manifold to which all the fibres of $g$ are now known to 
be isomorphic. Since $g$ is projective, we may choose a relatively very ample 
divisor $\Cal H$ on $\Cal Y$. Then, we may consider a morphism $\pi : \Cal B 
\rightarrow \Cal M$ defined by $b \mapsto [(\Cal Y_{b}, \Cal H \vert 
\Cal Y_{b})]$, where $\Cal M$ is an irreducible component of the coarse 
moduli scheme of polarised hyperk\"ahler manifolds. The existence of $\Cal M$ 
is guaranteed by Viehweg ([Vi]) as a very special case of his general Theory. 
Since $c_{1} : \text{Pic}(F) \rightarrow H^{2}(F, \Bbb Z)$ is injective for a 
hyperk\"ahler manifold, and since $H^{2}(F, \Bbb Z)$ is countable as a set, 
we see that $F$ admits only countably many polarisations. In particular, the 
isomorphism classes of $(\Cal Y_{b}, \Cal H \vert \Cal Y_{b})$ consists of at 
most countably many points. Therefore $\pi$ must be constant, because 
$\Cal B$ is irreducible and we are working over the uncountable base field 
$\Bbb C$. Set $m := \pi(\Cal B)$. Then $m$ is a single closed point of 
$\Cal M$. If $\Cal M$ happens to be a fine muduli 
space, then $g$ is obtained by the pull back of the universal family 
$u : \Cal U \rightarrow \Cal M$ by $\pi$. Hence, $g : \Cal Y \rightarrow 
\Cal B$ is the pull back of the single element $\Cal U_{m} \rightarrow \{m\}$ 
and therefore is globally trivial. Even in the general case, it is again 
known by 
E. Viehweg that we may work as if there were a universal family on $\Cal M$ 
in the following sense: 
There exist three objects: a reduced normal scheme $\Cal M'$; a finite 
automorphism group $\Gamma$ of $\Cal M'$ such that $\Cal M'/\Gamma = \Cal N$, 
where $\Cal N$ is the normalisation of the reduction $(\Cal M)_{\text{red}}$; 
and the family $u : \Cal U \rightarrow \Cal M'$ which is universal over the 
original $\Cal M$. (Refer to [Vi, Section 9, Theorem 9.25] for details. 
We just notice here that the assumption made in 
[Vi, Theorem 9.25] is satisfied if $\Cal M$ is a coarse moduli 
scheme of polarised manifolds of nef canonical classes, which the author 
learned from E. Viehweg.) Let us apply to our setting. Since $\Cal B$ is 
normal, $\pi$ factors through 
$\Cal N$. we denote this factorisation by $\pi_{\Cal N} : \Cal B \rightarrow 
\Cal N$. Since the natural map 
$\Cal N \rightarrow \Cal M$ is finite and since $\Cal B$ is irreducible, 
$\pi_{\Cal N}(\Cal B)$ is also a single closed point. Let us take an 
irreducible 
component $\Cal B'$ of the fibre product $\Cal B \times_{\Cal N} \Cal M'$ 
which dominates $\Cal B$. Then, the image of the fibre product map $\Cal B' 
\rightarrow \Cal M'$ is again a single closed point, because $\Cal M' 
\rightarrow \Cal N = \Cal M'/\Gamma$ is finite and $\Cal B'$ is irreducible. 
Let us denote this point by $m'$. Then the family $u_{\Cal B'} : 
\Cal U_{\Cal B'} \rightarrow \Cal B'$ obtained from 
$u : \Cal U \rightarrow \Cal M'$ by the fibre product map 
$\Cal B' \rightarrow \Cal M'$ is nothing but the pull 
back of the single $\Cal U_{m'} \rightarrow \{m'\}$. In particular, this is 
the trivial family: $p_{1} : \Cal B' \times F \rightarrow 
\Cal B'$. Since, $u : \Cal U \rightarrow \Cal M'$ is universal over $\Cal M$, 
the other pull back family $\Cal B' \times_{\Cal B} \Cal Y \rightarrow 
\Cal B'$ must coincide with the previous family and therefore is trivial. 
This already shows the isotriviality in the usual sense. Indeed, 
$\Cal B' \rightarrow \Cal B$ is finite and is even Galois, because so is 
$\Cal M' \rightarrow \Cal M'/\Gamma$. However, 
for our stronger version, we need to find a trivialisation by {\it \'etale} 
base change, which we will do from now by using the special fact that our 
family is a hyperk\"ahler family (but not much). Let us denote the Galois 
group of $\Cal B' \rightarrow \Cal B$ by $K$. Then, $g : \Cal Y \rightarrow 
\Cal B$ is the quotient $(p_{1} : \Cal B' \times F \rightarrow \Cal B')/K$, 
where the action of $K$ on $\Cal B' \times F$ is the action induced by the 
pull back $\Cal B' 
\rightarrow \Cal B$. Therefore, for each $k \in K$, the action of $k$ on 
$\Cal B' \times F$ is of the form 
$k : (b', f') \mapsto (k_{1}(b'), 
k_{2}(b', f'))$. Then, by fixing $k$ and varying $b' \in \Cal B'$, we obtain 
a  morphism $\Cal B' \rightarrow \text{Aut}(F)$ 
defined by $b' \mapsto k_{2}(b', *)$. Since $F$ is hyperk\"ahler, $\omega_{F}$ 
gives an isomorphism $T_{F} \simeq \Omega_{F}^{1}$. Therefore 
$H^{0}(T_{F}) \simeq H^{0}(\Omega_{F}^{1}) = 0$, where the last equality is 
because $\pi_{1}(F) = 0$. 
Hence $\text{Aut}(F)$ is discrete. Since $\Cal B'$ is irreducible, $\Cal B' 
\rightarrow \text{Aut}(F)$ is then constant. This means that the action of $k 
\in K$ is of the form $k : (b', f') \mapsto (k_{1}(b'), k_{2}(f'))$. 
In particular, we may speak of the homomorphism $p_{2} : K \rightarrow 
\text{Aut}(F)$ given by $k = (k_{1}, k_{2}) \mapsto k_{2}$. Assume 
that there exists a point $b' \in \Cal B'$ such that its stabiliser is 
non-trivial. Write this stabliser group by $H$. Then, $H$ acts on $F = \{b'\} 
\times F$ and induces an isomorphism $F \simeq F/H$, because 
$g : \Cal Y \rightarrow \Cal B$ 
is also a family of $F$. Since the caninical class of $F$ is trivial, the 
ramification formula implies that there are no ramification divisors 
of $F \rightarrow F/H$. Therefore, $F \rightarrow F/H$ is 
\'etale by the purity of the branch loci, because both source and target 
are smooth. Then, the quotient map 
$F \rightarrow F/H$ is an isomorphism, because $F/H \simeq F$ is simply 
connected. Hence, $\text{Im}(p_{2}\vert H : H 
\rightarrow \text{Aut}(F)) = \{id\}$. Therefore,  
the stabiliser group of each point of $\Cal B'$ is all contained in 
$M := \text{Ker}(p_{2} : K \rightarrow \text{Aut}(F))$. Note that $M$ is 
a normal subgroup of $K$. 
Now, by the definition of $M$, we see that $g : \Cal Y \rightarrow \Cal B$ is 
isomorphic to the quotient of $p_{1} : (\Cal B')/M \times F \rightarrow 
(\Cal B')/M$ by $K/M$. Moreover, by the construction, we see that 
$p_{1} : (\Cal B')/M \times F \rightarrow (\Cal B')/M$ has no point 
which admits a non-trivial stabiliser in $K/M$ any more. Therefore, the 
action of $K/M$ on $ (\Cal B')/M$ is free. 
Hence $p_{1} : (\Cal B')/M \times F 
\rightarrow (\Cal B')/M$ provides a desired trivialisation. \qed 
\enddemo
\demo{Proof of Corollary 3} Let us denote by $E_{w}$ the elliptic curve of 
period $w \in \Bbb H$. Choose a small neighbourhood $\varphi(\tau) \in 
\Delta_{2} \subset \Bbb H$ on which we have a family of elliptic curves 
$g : \Cal G \rightarrow \Delta_{2}$ such that $\Cal G_{w} = E_{w}$ and 
that $g$ has a section. Since $\varphi$ is holomorphic, we may also 
choose a small neighbourhood $\tau \in \Delta \subset \Bbb H$ such that 
$\varphi(\Delta) \subset \Delta_{2}$ and that there exists a family of 
elliptic curves $h : \Cal H \rightarrow \Delta$  
which satisfies that $\Cal H_{z} = E_{z}$ and admits a section. 
By pulling back $g : \Cal G \rightarrow \Delta_{2}$ by $\varphi$ 
and taking the fibre product, we obtain a family of abelian surfaces 
$a : \Cal A := \Cal H \times_{\Delta} \varphi^{*}\Cal G \rightarrow \Delta$. 
Here we have $\Cal A_{z} = E_{z} \times E_{\varphi(z)}$. Since $a$ admits 
a section, regarding this section as its $0$-section, we may consider the 
inversion $\iota$ of $\Cal A$ over $\Delta$. Dividing $\Cal A$ 
by $\iota$ and taking its crepant resolution, we obtain a family 
of Kummer surfaces $f : \Cal X \rightarrow \Delta$. By construction, 
we have $\Cal X_{z} = \text{Km}(E_{z} \times E_{\varphi(z)})$. This family 
$f$ is not a trivial family. Indeed, if $f$ is a trivial family, 
then there exists a K3 surface $S$ such that $S \simeq  \text{Km}(E_{z} 
\times E_{\varphi(z)})$ for all $z \in \Delta$. Note that the Kummer surface 
structures on $S$, that is, the isomorphism classes of abelian surfaces $A$ 
such that $S \simeq \text{Km}(A)$, are determined by the choices of $16$ 
disjoint smooth rational curves on $S$ and that there are at most countably 
many smooth rational curves on a K3 surface. Therefore, there exist at most 
countably many isomorphism classes of such $A$. Denote all of them by $A_{i}$ 
($i \in \Bbb N$) and set $A_{i} = \Bbb C^{2}/\Lambda_{i}$. For each $A_{i}$, 
the product structures 
on $A_{i}$, that is, the structures of decompositions $A_{i} = E_{i} \times 
F_{i}$, are also countably many, because subtori of $A_{i}$ are 
determined by the choices of sublattices of $\Lambda_{i}$. In conclusion, 
there are at most countably many isomorphism classes of pairs 
$(E, F)$ such that $S \simeq \text{Km}(E \times F)$. However, the set of 
the isomorphism classes of $E_{z}$ $(z \in \Delta)$ are uncountable, a 
contradiction. Therefore, our family $f : \Cal X \rightarrow \Delta$ 
is not trivial. Recall by [SM] that  $\rho(\Cal X_{z}) = 18$ if 
$E_{z}$ and $E_{\varphi(z)}$ are not isogenous and that 
$\rho(\Cal X_{z}) \geq 19$ if $E_{z}$ and $E_{\varphi(z)}$ are isogenous. 
Then, by Theorem 1, there exists a dense subset $\Cal S \subset \Delta$ 
such that $\rho(\Cal X_{s}) \geq 19$ for $s \in \Cal S$, that is, 
$E_{s}$ and $E_{\varphi(s)}$ are isogenous for $s \in \Cal S$. 
Therefore, any sequence in $\Cal S - \{\tau\}$ converging to $\tau$ 
satisfies our requirement. \qed 
\enddemo  

\demo{Proof of Corollary 6} Since both assumption and conclusion are 
compatible with shrinking of 
$\Delta$, we do this freely whenever it is convenient. Since $\pi : \Cal W 
\rightarrow \Delta$ is $\Bbb P^{1}$-bundle, $\pi$ is in particular a smooth 
morphism. Therefore, we may find three mutually disjoint analytic local 
sections of $\pi$ around $0$. Then $\Cal W = \Bbb P^{1} \times \Delta$ by the 
same argument as in the Harthshorne's book. By Theorem 1, there exists a 
dense countable subset $\Cal S \subset \Delta$ such that $\rho := 
\rho(\Cal X_{t})$ is constant for $t \in \Delta - \Cal S$ but 
$\rho(\Cal X_{s}) > \rho$ for $s \in \Cal S$. Let $\Cal D \subset \Bbb P^{1} 
\times \Delta$ be the discriminant locus of $\varphi$. In order to 
describe $\Cal D$, let us consider the Weierstrass model of $\varphi : \Cal X 
\rightarrow \Cal W$ with respect to $\Cal O$ and write the equation as 
$y^{2} = x^{3} + a(w, t)x + b(w, t)$, where $w$ is the inhomogeneous 
coordinate of $\Bbb P^{1}$ and $t$ is the coordinate of $\Delta$. Then 
$\Cal D$ is given by (the reduction of) the equation $4 a(w, t)^{3} + 27 
b(w, t)^{2} = 0 - (*)$. By construction, both $a(w, t)$ and $b(w, t)$ 
are polynomials with respect to $w$. Therefore the restriction map 
$\pi \vert \Cal D : \Cal D 
\rightarrow \Delta; (w, t) \mapsto t$ has at most finitely many such bad 
points $P \in \Cal D$ that $\pi \vert \Cal D$ is not smooth at $P$. Denote by 
$\Cal T \subset \Delta$ the set of the image of these bad points. This is 
then a finite set. In addition, since the type of non-multiple singular 
fibres are 
uniquely determined by the local monodromy, the singular fibres of the 
fibrations $\varphi_{t} : \Cal X_{t} \rightarrow \Bbb P^{1}$ ($t \in \Delta - 
\Cal T$) are exactly the same regardless of $t$. Write them by $T_{i}$ 
($i = 1, \cdots n$) and denote by $m_{i}$ the number of the irreducible 
components of $T_{i}$. Then, by Shioda's formula [Sh1], we have 
$r(\varphi_{t}) = \rho(\Cal X_{t}) - 2 - \sum_{i = 1}^{n} (m_{i} - 1)$ 
for $t \in \Delta - \Cal T$. 
Set $\Cal S' := \Cal S - \Cal T$. Then, at each point of $\Cal S'$, the 
function $r(t)$ is strictly upper semi-continuous in the $I$-topology. 
Since $\Cal S$ is countable and dense (in the classical topology) 
and $\Cal T$ is a finite set, $\Cal S'$ is also countable 
and dense. Therefore, this $\Cal S'$ provides a desired set. \qed 
\enddemo
\demo{Proof of Corollary 8} First we shall show the exsistence of a sequence 
in the statement. We may assume that $r \leq 17$. 
Let us consider the Kuranishi family 
$k : (J \subset \Cal U) \rightarrow (0 \in \Cal K)$ of $J$. This is a germ of 
the universal deformation of $J$ and is known to be smooth of dimension 20. 
Threfore, $\Cal K$ is realised as an open neighbourhood of $0 \in 
H^{1}(J, T_{J})$ and is assumed to be a small polydisk in $\Bbb C^{20}$. 
Then $R^{2}k_{*} \Bbb Z_{\Cal U}$ is a constant system on $\Cal K$. Let us 
fix a marking $\tau 
: R^{2}k_{*} \Bbb Z_{\Cal U} \simeq \Lambda \times \Cal K$, where 
$\Lambda = \Lambda_{\text{K3}}$, and consider as before the resulting period 
map 
$$p : \Cal K \rightarrow \Cal D = \{ [\omega] \in \Bbb P(\Lambda \otimes 
\Bbb C) \vert (\omega.\omega) = 0, (\omega.\overline{\omega}) > 0 \} 
\subset  \Bbb P(\Lambda \otimes \Bbb C) \simeq \Bbb P^{21}.$$ 
Then, $p$ is a local isomorphism by the local Torelli Theorem and by the 
fact that $\text{dim}\Cal K = \text{dim} \Cal D (= 20)$. Therefore, 
we may identify $\Cal K$ with an open neighbourhood $\Cal U \subset \Cal D$ 
of $p(0)$ by $p$. Since our argument is completely local, 
by abuse of notation, we 
write this $\Cal U$ again by $\Cal D$ and identify therefore $\Cal K = \Cal D$ 
by $p$. Let us write a general fiber of $\varphi : J \rightarrow 
\Bbb P^{1}$ by $E$ and choose an integral basis $S_{i}$ ($i = 1, \cdots, r$) 
of $MW^{0}(\varphi)$, where $r := r(\varphi)$. Then $S_{i}$ and the zero 
section $O$ are all non-singular rational curves and $E$ is an elliptic curve 
such that $(S_{i}.E) = (O.E) = 1$.  For this purpose, we may 
assume that 
$r \leq 17$. By the definition of the Mordell-Weil lattice 
$(MW(\varphi), \langle *, * \rangle)$, we have the (minus sign of) 
isometric injective homomorphism  $\iota : MW^{0}(\varphi) \hookrightarrow 
NS(J)$ given by $S \mapsto S - O$ [Sh2]. Then, $\langle S_{i}, S_{i} 
\rangle = 4 + 2(S_{i}.O)$. Note also that $E$, $O$, $S_{i}$ are linearly 
independent in $NS(J)$ [Sh1]. Let us consider $(r+2)$ elements in 
$\Lambda_{\text{K3}}$ given by $e := \tau([E])$, $o := \tau([O])$ and 
$s_{i} := \tau([S_{i}])$. Then, these are also linearly indepedent in 
$\Lambda$. Consider the subset $\Cal L$ of $\Cal D$ defined by 
$(e.*) = (o.*) = (s_{i}.*) = 0$. This is a smooth analytic subset of 
$\Cal D$ of dimension $20 - (r + 2)$ which contains $p(0)$. 
This follows from the non-degeneracy of $(*.*)\vert 
\langle o, e, s_{i} \rangle^{\perp}$ and the Lefschetz $(1,1)$-Theorem. 
Through the identification made 
above, we may regard $0 \in \Cal L \subset \Cal K$. Then we may speak of the 
family $\tilde{j} : \tilde{\Cal J} \rightarrow \Cal L$ obtained as the 
restriction of 
$k : \Cal U \rightarrow \Cal K$ to $\Cal L$. Then by [Ko, Theorem 14] 
(see also [MM]) or by [Hu1, Section 1 (1.14)] in more sophisticated 
terminology, we see that 
$\Cal L$ is the locus on which the invertible sheaves $\Cal O_{J}(E)$, 
$\Cal O_{J}(O)$, $\Cal O_{J}(S_{i})$ on $J$ lift to invertible sheaves 
$\Cal E$, $\Cal O$ and $\Cal S_{i}$ on the whole space $\tilde{\Cal J}$. 
Since $20 - (r+2) \geq 1$ by $r \leq 17$, we can take a 
sufficiently small disk $0 \in \Delta \subset \Cal L$. Then we may consider 
the induced family $j : \Cal J \rightarrow \Delta$. We denote the restrictions 
of $\Cal E$, $\Cal O$ and $\Cal S_{i}$ on $\Cal J$ by the same letter. We also 
shrink $\Delta$ freely whenever it is convenient. 
Note that $\chi(\Cal O_{J}(S_{i})) = 1$, $h^{0}(\Cal O_{J}(S_{i})) = 1$ 
and $h^{q}(\Cal O_{J}(S_{i})) = 0$ for $q > 0$, because $S_{i}$ is a smooth 
rational curve on a K3 surface. Then by applying the upper 
semi-continuity of coherent sheaves and by the Theorem of cohomology, we 
see that $j_{*}\Cal S_{i}$ are invertible sheaves which satisfy the base 
change property. Then $(j_{*}\Cal S_{i}) \otimes \Bbb C(0) 
\simeq H^{0}(\Cal O_{J}(S_{i}))$. Therefore, by Nakayama's Lemma, 
all of $C_{i}$ lift not only as invertible sheaves but also as effective 
divisors on $\Cal J$. By abuse of notation, we denote these divisors again 
by $\Cal S_{i}$. Since the smoothness is an open condition for a proper 
morphism, 
$\pi \vert \Cal S_{i} : \Cal S_{i} \rightarrow \Delta$ is also smooth. 
Combining this with the fact that small 
deformation of $\Bbb P^{1}$ is again $\Bbb P^{1}$, we see that $S_{i, t} := 
\Cal S_{i} \vert \Cal J_{t}$ is again a smooth rational curve on $\Cal J_{t}$ 
for all $t \in \Delta$. The same holds for $\Cal O_{t} := \Cal O \vert 
\Cal J_{t}$. Note that $\chi(\Cal O_{J}(E)) = 2$, $h^{0}(\Cal O_{J}(E)) = 2$ 
and $h^{q}(\Cal O_{J}(E)) = 0$ for $q > 0$, because $E$ is an elliptic curve 
on a K3 surface. Then, $h^{q}(\Cal E \vert \Cal J_{t}) = 0$ and 
$h^{0}(\Cal E \vert \Cal J_{t}) = 2$. Therefore, $j_{*} \Cal E$ is a locally 
free sheaf of rank 2 which satisfies the base change property. In 
particular, $j^{*}j_{*}\Cal E \vert J = H^{0}(\Cal O_{J}(E))$. Since 
$\Cal O_{J}(E)$ is globally generated, we see again by Nakayama's Lemma 
that the natural map $j^{*}j_{*}\Cal E \rightarrow \Cal E$ is also surjective. 
Therefore we may 
take a morphism $\Phi : \Cal J \rightarrow \Cal W$ over $\Delta$ associated 
to this surjection. Then, by the base change property, we find that the 
restriction $\Phi_{t} : \Cal J_{t} \rightarrow 
\Cal W_{t}$ coincides with the morphism given by the surjection 
$H^{0}(\Cal E \vert \Cal J_{t}) \otimes \Cal O_{\Cal J_{t}} \rightarrow 
\Cal E \vert \Cal J_{t}$. This is an elliptic fibration by 
$h^{0}(\Cal E \vert \Cal J_{t}) = 2$ and by the adjunction formula on a 
K3 surface. Therefore, the factorisation $\Phi : \Cal J \rightarrow \Cal W$ 
makes $j : \Cal J \rightarrow \Delta$ a family of elliptic surfaces over 
$\Delta$. By the invariance of the intersection number, we have 
$(\Cal S_{i, t}. \Cal E_{t}) = (S_{i}.E) = 1$. Therefore, $S_{i,t}$ is also 
a section of $\Phi_{t}$. The same holds for $\Cal O_{t}$. Therefore 
$\Phi : \Cal J \rightarrow \Cal W$ makes $j : \Cal J \rightarrow \Delta$ 
a family of Jacobian K3 surfaces with zero section $\Cal O$. Moreover, 
by passing to 
the Weierstrass family over $\Delta$ given by $\Cal O$ and using the 
characterisation of $MW^{0}(\varphi)$ that $S \in MW(\varphi)$ is in 
$MW^{0}(\varphi)$ if and only if $S$ does not meet the singular points of the 
Weierstrass model, we also see that $\Cal S_{i, t}$ are all in 
$MW^{0}(\Phi_{t})$. In addition, the 
intersection matrix of $\Cal E_{t}$, $\Cal O_{t}$, $\Cal S_{i, t}$ are the 
same as the one for $E$, $O$, $S_{i}$ in $\Lambda$ and is then hyperbolic. 
Therefore, $\Cal E_{t}$, $\Cal O_{t}$, $\Cal S_{i, t}$ are also linearly 
independent in $H^{2}(\Cal J_{t}, \Bbb Z)$. Hence so are in $NS(\Cal J_{t})$. 
Thus by the 
injection $MW^{0}(\Phi_{t})  \hookrightarrow NS(\Cal J_{t})$ quoted above, we 
see that $S_{i,t}$ are also linearly independent in $MW^{0}(\Phi_{t})$. In 
particular, $r(\Phi_{t}) \geq r$ for all $t \in \Delta$. Since the base space 
$\Delta$ is chosen in the Kuranishi space, our family $j : \Cal J \rightarrow 
\Delta$ is not trivial. Therefore, by Corollary 6, there exists $t_{0} \in 
\Delta$ such that $r(t)$ is strictly upper semi-continuous at $t_{0}$. In 
particular,  $r(t_{0}) > r$. By the invariance intersection and 
by the relation between $\langle *, * \rangle$ and $(*, *)$ quoted above, 
we see that the map $a : MW^{0}(\varphi) 
\rightarrow MW^{0}(\Phi_{t_{0}})$ given by $S_{i} \mapsto \Cal S_{i, t_{0}}$ 
is then an isometric injection. If $r(t_{0}) = r + 1$, then we may define 
$\varphi_{r+1} : J_{r+1} \rightarrow \Bbb P^{1}$ to be this Jacobian K3 
surface $\Phi_{t_{0}} : \Cal J_{t_{0}} 
\rightarrow \Bbb P^{1}$. Let us treat the case where 
$r(t_{0}) \geq r + 2$. Since $18 \geq r(t_{0})$, we have $16 \geq r$. 
For simplicity, 
we abbreviate $\Phi_{t_{0}} : \Cal J_{t_{0}} 
\rightarrow \Bbb P^{1}$ and $r(t_{0})$ by 
$\varphi' : J' \rightarrow \Bbb P^{1}$ and $r'$ respectively. We denote the 
image of the basis $S_{i}$ ($1 \leq i \leq r$) of $MW^{0}(\varphi)$ 
in $MW^{0}(\varphi')$ by the same letters $S_{i}$ and take 
$T_{j} \in MW^{0}(\varphi')$ $j = r+1, r+2, \cdots, r'$ such that 
$S_{i}$ and $T_{j}$ form a basis of $MW^{0}(\varphi') \otimes 
\Bbb Q$ over $\Bbb Q$. (Here note that our embedding might not be primitive 
so that we can not prolong $S_{i}$ to an integral basis of 
$MW^{0}(\varphi')$ in general.) Let us consider the Kuranishi space 
$\Cal K'$ of $J'$ and take the subspace $\Cal L' \subset \Cal K'$ 
defined by the fibre 
class $E'$ of $\varphi'$, the zero section $O$, all of $S_{i}$, and 
$T_{r+1}$, and denote by $j' : \Cal J' \rightarrow \Cal L'$ the family 
induced by the Kuranishi family as before. Then, 
$\text{dim}\Cal L' = 20 - (2 + r + 1) > 0$, because $E'$, $O$, 
$S_{i}$ and $T_{r+1}$ are linearly independent in $H^{2}(J', \Bbb Z)$ 
and because $r \leq 16$. 
In addition, considering $\Cal L'$ as a 
subspace in the period domain under the identification made as before, and 
applying the same argument as in (1.1) based on the holomorphicity of the 
preiod mapping and the Lefschetz (1,1)-Theorem, we see that the 
N\'eron-Severi group of any general fibre $\Cal J_{t}'$ in the sense of the 
$I$-topology (Here we define a non-empty open set as a complement of 
countably many proper analytic subsets) is isomorphic to the primitive 
closure of $\Bbb Z \langle E_{18}, O, S_{i}, T_{r+1} \rangle$ in 
$H^{2}(J', \Bbb Z)$. In particular, $\rho(\Cal J_{t}') = r + 3$ 
for general $t$. Moreover, 
by the same argument as above, we find that this family becomes a 
family of Jacobian K3 surfaces 
$\Cal J' \overset \Phi' \to \longrightarrow \Cal W' \rightarrow \Cal L'$ 
such that each fibre $\Phi_{t}' : \Cal J_{t}' \rightarrow \Cal W_{t}'$ 
satisfies that $MW^{0}(\varphi) 
\subset MW^{0}(\Phi_{t}')$ and that
$r(\Phi_{t}') \geq r + 1$. On the other hand, we have $r(\Phi_{t}') \leq 
r + 1$ for general $t$ by Shioda's formula and by $\rho(\Cal J_{t}') = r + 3$. 
Then we have $r(\Phi_{t}') = r + 1$ and may define $\varphi_{r+1} : J_{r +1} 
\rightarrow \Bbb P^{1}$ to be $\Phi_{t}' : \Cal J_{t}' \rightarrow 
\Cal W_{t}'$ for general $t$. 
The first statement of Corollary 8 now follows from induction on $(18 - r)$. 
Next we shall show the middle statement of Corollary 8. Let $\phi' : S' 
\rightarrow \Bbb P^{1}$ be a Jacobian K3 
surface such that $\text{rank}(\phi') = 18$ and $M$ a sublattice of 
$MW(\phi')$. Then, by taking a general point of the locus of the Kuranishi 
space defined by the basis of $M$, zero section of 
$\phi'$ and general fibre of $\phi'$ as before, we obtain a Jacobian K3 
surface 
$\phi : S \rightarrow \Bbb P^{1}$ such that $M \subset MW^{0}(\phi)$ 
and $r(\phi) = r$. Finally, we show the last assertion of Corollary 8. 
Assume that $M \subset MW^{0}(\phi)$ and is of finite index. 
Since the pairing $\langle * , * \rangle$ is integral valued on 
$MW^{0}(\phi)$ [Sh2], we have then $M \subset MW^{0}(\phi) \subset M^{*}$. 
Since $M \subset M^{*}$ is of finite index, the possibilities of 
$MW^{0}(\phi)$ is then only finitely many. By [Sh2], we have also 
$MW^{0}(\phi) \subset MW(\phi)/\text{(torsion)} \subset 
MW^{0}(\phi)^{*}$. Therefore each given $MW^{0}(\phi)$
also recovers $MW(\phi)/\text{(torsion)}$ up to finitely many ambiguities. 
Now it is sufficient to show the boundedness of the torsion subgroups of 
Jacobian K3 surfaces. For those which have non-constant $j$-invariant, 
the result 
follows from the classification due to Cox [Co]. Let us consider the case 
where the $j$-invariant is constant. Note that a Jacobian 
K3 surface always admits at least one singular fibre, because 
its topological Euler number is positive. Therefore, by the classification of 
the singular fibres whose $j$-values are not $\infty$ and 
by the general fact that 
the specialisation map $MW(\phi)_{\text{torsion}} \rightarrow 
(\phi^{-1}(t))_{\text{reg}}$ is 
injective, we see that the possible torsion groups are: 
$0$, $\Bbb Z/2$, $\Bbb Z/3$, $\Bbb Z/4$ or $(\Bbb Z/2)^{\oplus 2}$. 
This implies the result. 
\qed 
\enddemo
\par 
\vskip 4pt 
\head
{\S 3. Construction of Examples}
\endhead 
We shall give explicit construction of Examples in Introduction. 
\demo{Construction of Example 4} By [OZ] based on [SI], there exists a K3 
surface $S$ of $\rho(S) = 20$ which contains 19 smooth rational curves $C_{i}$ 
($i = 1, 2, \cdots, 19$) whose configuration is of 
Type $D_{19}$ and is primitive in $NS(S)$. Consider the Kuranishi space of 
$S$ and take the locus defined by $(C_{i}.*) = 0$ ($1 \leq i \leq 19$) as 
in the proof of Corollary 8. Then this is a smooth curve and we may take 
a small curve 
$\Delta$ inside this curve and speak of the family $f : \Cal X \rightarrow 
\Delta$ induced by the Kuranishi family. Then, by the same argument as 
in Corollary 8, we see that $\rho(t) \geq 19$ for $t \in \Delta$ and that if 
$\rho(t) = 19$, then $NS(\Cal X_{t}) \simeq D_{19}$. In the last case, we 
have $a(\Cal X_{t}) = 0$, because $D_{19}$ is negative definite. 
Note that $f$ is not trivial. Then by Theorem 1, we have the dense 
countable set $\Cal S$ such that $\rho(\Cal X_{s}) = 20$ for $s \in \Cal S$ 
and $\rho(\Cal X_{t}) = 19$ for $t \not\in \Cal S$. Since a K3 surface of the 
maximal Picard number 20 is algebraic, this family satisfies all the 
properties required. \qed 
\enddemo 
\demo{Construction of Example 5} Applying the construction 
in Corollary 3 for $\varphi(z) \equiv \sqrt{-1}$, we obtain a family of 
Kummer surfaces $f : \Cal X \rightarrow \Delta$ such that 
$\Cal X_{z} = \text{Km}(E_{z} \times E_{\sqrt{-1}})$. Then by [SM], we have 
$\rho(\Cal X_{z}) = 18$ for 
$z \not\in \Bbb Q(\sqrt{-1})$, and 
$\rho(\Cal X_{z}) = 20$ for $z \in \Bbb Q(\sqrt{-1})$. Since 
$\Bbb Q(\sqrt{-1})$ is dense in $\Delta$, this provides an example we seeked. 
\qed 
\enddemo
\demo{Construction of Example 7} Let us first consider the family of rational 
Jacobian surfaces with rational double points $h : \Cal Z \rightarrow 
\Bbb P^{1} \times \Delta_{u} \rightarrow \Delta_{u}$ defined by the 
Weierstrass equation $y^{2} = x^{3} + ux + s^{5}$. Here $u$ is the 
coordinate of $\Delta$ and $s$ is the inhomogeneous coordinate of 
$\Bbb P^{1}$. Then either by the N\'eron algorithm or by a direct 
calculation, we can easily check the following fact: 
$\Cal Z$ is smooth; $\Cal Z_{u}$ ($u \not= 0$) is smooth and 
$\Cal Z_{u} \rightarrow \Bbb P^{1}$ has singular fibres of Type $I_{1}$ over 
$4u^{3} + 27s^{10} = 0$ and of Type $II$ over $s = 
\infty$; and $\Cal Z_{0} \rightarrow \Bbb P^{1}$ has one singular point of 
type $E_{8}$ over $s = 0$ 
and has a singular fibre of Type $II$ over $s = \infty$. Therefore taking an 
appropriate finite covering $\Delta_{v} \rightarrow \Delta_{u}$ and a 
simultaneous resolution of the pull back family, we obtain a family of smooth 
rational Jacobian surfaces $g : \Cal Y \rightarrow \Bbb P^{1} \times 
\Delta_{v} \rightarrow \Delta_{v}$ such that $\Cal Y_{v} \rightarrow 
\Bbb P^{1}$ ($v \not = 0$) has 10 singular fibre of Type $I_{1}$ and one 
singular fibre of Type $II$, and $\Cal Y_{0} \rightarrow \Bbb P^{1}$ has 
one singular fibre of Type $II^{*}$ and one singular fibre of Type $II$. Then 
by Shioda's formula, we have $r(v) = 8$ for $v \not= 0$ and $r(0) = 0$. 
Let us choose large number $M$ such that the divisor $s = M$ on 
$\Bbb P^{1} \times \Delta_{v}$ does not meet the discriminant locus. 
This is possible by the description above. Let us take the double covering 
$\Bbb P^{1} \times \Delta_{t} \rightarrow \Bbb P^{1} \times \Delta_{v}$ 
ramified over $s = M$ and $s = \infty$ and consider the relatively 
minimal model $f : \Cal X \rightarrow \Bbb P^{1} \times \Delta_{t} 
\rightarrow \Delta_{t}$ of the pull back family 
$f' : \Cal X' \rightarrow \Bbb P^{1} \times \Delta_{t} \rightarrow \Delta_{t}$ 
over $\Delta_{t}$. Note that $\Cal X'$ is equi-singular along the preimage of 
the cuspidal points of fibres $(\Cal Y_{v})_{\infty}$ 
of $\Cal Y_{v} \rightarrow \Bbb P^{1}$. Then, by the monodromy calculation, 
we see that $f : \Cal X \rightarrow \Bbb P^{1} \times \Delta_{t} 
\rightarrow \Delta_{t}$ is a smooth family of Jacobian K3 surfaces such that 
$\Cal X_{0} \rightarrow \Bbb P^{1}$ has two singuler fibres of Type $II^{*}$ 
and one singular fibre of Type $IV$; and $\Cal X_{t} \rightarrow \Bbb P^{1}$ 
($t \not= 0$) has 20 singular fibres of Type $I_{1}$ and one singular 
fibre of Type $IV$. Note also that $r(t) \geq r(v) = 8$ for $t \not= 0$. 
On the 
other hand, again by Shioda's formula, we have $20 \geq \rho(t=0) = 2 + 
r(t=0) + (9 -1) + (9 -1) + (3 -1)$. Therefore $\rho(0) = 20$ and 
$r(0) = 0$ for the 
centre fibre of $f$. Moreover, this family $\Cal X \rightarrow \Delta$ is not 
trivial even as a family of K3 surfaces. Indeed, otherwise, we have $\Cal X 
\simeq \Cal X_{0} \times \Delta_{t}$. Since $\text{Pic}(\Cal X_{0})$ is 
discrete, it is impossible for $\Cal X_{0}$ to admit a family of elliptic 
fibrations which vary continuously. Then, our family 
$\Cal X \rightarrow \Delta$ must be also trivial 
as a family of elliptic fibre spaces if $\Cal X 
\simeq \Cal X_{0} \times \Delta_{t}$ . However, this contradicts the fact 
that the type of singular fibre of $\Cal X_{0}$ and $\Cal X_{t}$ are 
different. Therefore our family is not trivial even as a family of K3 
surfaces. Hence, by (1.2), we have $\rho(t) < 20$ 
for general $t$ in the $I$-topology. 
\qed 
\enddemo 
\remark{Correction of [OS]} Taking this place, we shall quote a correction of 
mistakes of [OS Main Theorem] pointed out by Y. Yamada: 
\par 
No 12 : the first line of the matrix of $E(K)^{0}$ should be $4, -1, 0, 1$.
\par 
No 32 : $E(K)^{0}$ should be $\left(\matrix 4 & -2 \\ -2 & 4 \endmatrix 
\right)$.
\par
No 70 : $E(K)$ should be $\Bbb Z/4\Bbb Z$. \qed
\endremark

\enddocument